\documentclass{conm-p-l}

\usepackage{amsmath}
\usepackage{amscd}
\usepackage{amssymb}

\newtheorem{theorem}{Theorem}[section]

\theoremstyle{definition}

\newtheorem{proposition}[theorem]{Proposition}
\newtheorem{remark}[theorem]{Remark}

\theoremstyle{remark}

\newcommand{\ale}{ \C^{2} // \G}
\newcommand{\alen}{ (\ale)^{[n]} }
\newcommand{\C}{{\mathbb C}}
\newcommand{\ch}{{\rm ch}}

\newcommand{\diag}{ {\rm diag}}
\newcommand{\End}{ {\rm End}}
\newcommand{\Id}{ {\rm Id}}
\newcommand{\Supp}{ {\rm Supp}}
\newcommand{\g}{\gamma}
\newcommand{\G}{\Gamma}
\newcommand{\Gm}{\G_m}

\newcommand{\Gn}{\G_n}
\newcommand{\hquiver}{Y_{\G,n}}
\newcommand{\hilq}{Y // \G}
\newcommand{\hilqgn}{ \C^{2n} // \Gn }

\newcommand{\Hn}{H^*(\Xn)}
\newcommand{\Hom}{{\rm Hom}}
\newcommand{\Hx}{\mathbb H_X}
\newcommand{\quotuniv}{{\mathcal U}_{\G, n} }
\newcommand{\RG}{\mathbb R_{\G}}
\newcommand{\Xm}{X^{[m]}}
\newcommand{\Xmn}{X^{[m+n]}}
\newcommand{\Xn}{X^{[n]}}
\newcommand{\vac}{|0\rangle}
\newcommand{\w}{\tilde}
\newcommand{\wt}{\xi}
\newcommand{\Z}{ {\mathbb Z} }

\numberwithin{equation}{section}

\begin{document}

\title
[Hilbert schemes and wreath products]{Algebraic structures behind
Hilbert schemes and wreath products}
\author{Weiqiang Wang}
\address{Department of Mathematics\\
North Carolina State University\\ Raleigh, NC 27695;\\
current address: Department of Mathematics\\
University of Virginia  \\Charlottesville, VA 22904.}
\email{ww9c@virginia.edu}
\thanks{The author is supported in part by an NSF Grant
and an FR\&PD grant from NCSU}

\subjclass{Primary 14C05, 20C05; Secondary 17B67, 17B68, 17B69}

\begin{abstract}
In this paper we review various strikingly parallel algebraic
structures behind Hilbert schemes of points on surfaces and
certain finite groups called the wreath products. We explain
connections among Hilbert schemes, wreath products,
infinite-dimensional Lie algebras, and vertex algebras. As an
application we further describe the cohomology ring structure of
the Hilbert schemes. We organize this paper around several
general principles which have been supported by various works on
these subjects.
\end{abstract}
\maketitle
\tableofcontents

\section{Introduction}

The purpose of this paper is to describe the analog and
connections between two seemingly unrelated subjects. The first
one concerns about a wonderful geometric object, namely, the
Hilbert scheme $\Xn$ of $n$ points on a (quasi-)projective surface
$X$. It has been well known \cite{Fog} that $\Xn$ is non-singular
of complex dimension $2n$. The Hilbert-Chow morphism from $\Xn$ to
the symmetric product $X^n/S_n$ is a semismall crepant resolution
of singularities. The second one is a finite group called the
wreath product $\Gn$, namely, the semi-direct product between the
symmetric group $S_n$ and the product group $\G^n$ of a finite
group $\G$. The representation theory of wreath products was first
developed by Specht \cite{Spe}, also cf. e.g. \cite{Mac, Zel}. For
$\G$ trivial, $\Gn$ reduces to $S_n$. Many nice features of the
representation theory of symmetric groups remain to hold in the
setup of wreath products.

It turns out that there exist deep connections among the geometry
of Hilbert schemes, representation theory of wreath products, and
vertex algebras. In recent years, these subjects have attracted
various people with diversified backgrounds such as algebraic
geometry, representation theory, combinatorics, and mathematical
physics. In this mostly expository paper, we attempt to organize our
discussions around four main points which have been supported by
various works on these subjects. (It also includes some original materials
from the unpublished paper \cite{Wa2} and further clarification).
It is our hope that many topics
under discussion may attain somewhat better understanding in this
way and these principles may serve again as a helpful guide in
uncovering new structures.

The {\em first point} is that one should study the Hilbert schemes
$\Xn$ (respectively the wreath products $\Gn$) for all $n \ge 0$
simultaneously. It was not too long ago when people started to
study the geometry of Hilbert schemes in such a way. G\"ottsche
\cite{Got1, Got2} calculated the Betti numbers for the Hilbert
scheme $\Xn$ associated to a surface $X$ and presented it in a
beautiful formula in terms of a generating function for all $n \ge
0 $ together. Motivated by G\"ottsche's formula, quiver varieties
\cite{Na1} and Vafa-Witten's paper on $S$-duality \cite{VW},
Nakajima \cite{Na2, Na3} constructed a Heisenberg algebra
associated to the lattice $H^*(X, \Z)/{tor}$ in terms of the
correspondence varieties which act irreducibly on the direct sum
$\Hx$ of the cohomology groups of $\Xn$ for all $n \ge 0$. Similar
results were obtained by Grojnowski \cite{Gro}. On the other hand,
it has been known for a long time that one should study the
representations of symmetric groups $S_n$, or more generally the
wreath products $\Gn$, for all $n \ge 0$ together, cf. \cite{Spe,
Zel, Mac}. A Heisenberg algebra associated to the lattice
$R_{\Z}(\G)$ was constructed by the author \cite{Wa1} (also cf.
\cite{FJW1}) acting irreducibly on the direct sum $\RG$ of
representation rings of $\Gn$ for all $n \ge 0$. Here $R_{\Z}(\G)$
denotes the integral span of the irreducible characters of $\G$.

The {\em second point} is that the geometry/invariants of a
suitable resolution of singularities of an orbifold should be
compared to the equivariant counterparts of the orbifold. This was
largely stimulated by the study of the orbifold string theory
where the notion of orbifold Euler numbers was introduced
\cite{DHVW}, and it may also be viewed as a McKay correspondence
in a broad sense, cf. \cite{Rei1, Rei2} and references therein.
Given an orbifold, one may study the associated equivariant
$K$-group, equivariant derived category, orbifold cohomology,
orbifold Euler number, orbifold Hodge number, orbifold elliptic
genera etc, and compare them with their counterparts on a suitable
resolution, cf. \cite{DHVW, Zas, DMVV, Bat, BKR, Wa2, CR} and
references therein. The Euler numbers and Hodge numbers of Hilbert
schemes have been computed in \cite{Got1, GS}. The Hilbert-Chow
morphism $\Xn \rightarrow X^n/S_n$ provides a great example which
matches invariants on a suitable resolution with those on the
orbifold, cf. \cite{HH, Got3, DMVV, Zhou}.

The analog and connections between Hilbert schemes and wreath
products were pointed out by the author \cite{Wa1}. Given a space
$Y$ with an action by a finite group $\G$, the wreath product
$\Gn$ acts on the $n$-th direct product $Y^n$ in a canonical way.
Below we assume $\tau: X \rightarrow Y/\G$ is a resolution of
singularities of the orbifold $Y/\G$ and assume both $X$ and $Y$
are complex surfaces. We observe that there is a naturally induced
resolution of singularities $\tau_n: \Xn \rightarrow Y^n/\Gn$. The
{\em third point} is that whenever the resolution $\tau: X
\rightarrow Y/\G$ is `good' in a suitable sense then the
resolution $\tau_n: \Xn \rightarrow Y^n/\Gn$ is `good'. A variant
of this can be formulated as follows: whenever a certain
`reasonable' statement can be made relating $X$ and $Y/\G$ then
the corresponding statement relating $\Xn$ and $Y^n/\Gn$ should
hold.

For example, if a resolution $\tau: X \rightarrow Y/\G$ is
crepant, respectively semismall, then so is the corresponding
resolution $\tau_n: \Xn \rightarrow Y^n/\Gn$. If the Euler number
of $X$ is equal to the orbifold Euler number of the orbifold
$Y/\G$, then the Euler number of $\Xn$ is equal to the orbifold
Euler number of the orbifold $Y^n/\Gn$, cf. \cite{Wa1}. The
statement remains to be true by considering Hodge numbers instead
of Euler numbers, and it is conjectured the statement is also true
for the elliptic genera, cf. Wang-Zhou \cite{WaZ}. We remark that
Borisov and Libgober \cite{BL} have formulated mathematically the
notion of an orbifold elliptic genera and verified a conjectural
formula in \cite{WaZ} for the orbifold elliptic genera of
$Y^n/\Gn$. When $\G$ is trivial and $X$ equals $Y$, the third
principle above simply says the Hilbert-Chow resolution is `good'
in all aspects mentioned in the previous paragraphs. Since one can
easily construct various examples of good resolutions $\tau: X
\rightarrow Y/\G$, this principle provides in a tautological way
numerous new higher dimensional examples of good resolutions. To
our best knowledge, (except a few isolated cases) all the known
higher dimensional examples of good resolutions arise in such a
way.

This above discussion `explains' why it is natural to study the
direct sum of the equivariant $K$-groups $K_{\Gn}(Y^n) \bigotimes
\C$ and why this should be parallel to the study of $\Hx$
\cite{Wa1} (also cf. \cite{Seg2} and a footnote in \cite{Gro} for
the important case when $\G$ is trivial). To simplify the
discussion, we will refrain ourselves to discuss the
representation rings $R(\Gn)$ (instead of the equivariant $K$
groups) with connections to the cohomology groups of Hilbert
schemes. In an important case when $\G$ is a finite subgroup of
$SL_2(\C)$, $Y$ is the affine plane $\C^2$, and $X$ is the minimal
resolution of singularities of the simple singularity $\C^2/\G$,
considering only the representation rings $R(\Gn)$ does not really
lose information since the equivariant $K$ group $K_{\Gn}(\C^{2n})
\bigotimes \C$ is isomorphic to $R(\Gn)$ by the Thom isomorphism.
This lead the author to propose a group theoretic realization,
which can be viewed as a new variation of McKay correspondence and
have been subsequently developed jointly with I.~Frenkel and Jing
\cite{FJW1}, of the homogeneous vertex representation of the
affine Kac-Moody algebra of ADE types and its toroidal
counterpart (cf. \cite{FK, Seg1, MRY}). We remark that a weighted
bilinear form on $R(\Gn)$ introduced in \cite{FJW1} in a group
theoretic manner affords a natural interpretation in terms of the
Koszul-Thom complex \cite{Wa2}. This approach to the McKay
correspondence is parallel to a geometric realization in terms of
Hilbert schemes $\Xn$ when $X$ is the minimal resolution $\ale$ of
the simple singularity $\C^2/\G$, cf. \cite{Na3}.

The {\em fourth point} is that both the Hilbert schemes and wreath
products have deep connections with the theory of vertex algebras
\cite{Bor1}\footnote{Vertex algebras have played an important role
in many different fields, including string theory,
infinite-dimensional Lie algebras, the Monster group, the
moonshine conjecture, mirror symmetry, moduli space of algebraic
curves, and others, cf. \cite{Kac1, FLM, Bor2, Kac2, Fre} and the
references therein.}, and it is beneficial to study them in a
parallel way. The shift from a representation to its underlying
symmetry algebra provides new insights into various structures of
the representation. In the framework of Hilbert schemes, the first
indication of connections with vertex algebras comes from the
construction of the Heisenberg algebra (cf. \cite{Na2, Gro, Na3}).
The Heisenberg algebra construction was used by de Cataldo and
Migliorini \cite{dCM} to study the geometry of Hilbert schemes
from a novel viewpoint. Lehn \cite{Lehn} realized geometrically
the Virasoro algebra and applied it to study connections between
Heisenberg generators and cup products on the cohomology groups of
Hilbert schemes. In our joint work with Li and Qin \cite{LQW1,
LQW2}, we have further developed the connection between vertex
algebras and Hilbert schemes and used it to obtain new results on
the cohomology ring structure of $\Xn$ associated to an arbitrary
projective surface $X$ (which in general has been unaccessible by
classical methods in algebraic geometry). On the other hand, in
the framework of wreath products, there has been a group theoretic
construction of the Virasoro algebra given by I.~Frenkel and the
author \cite{FW} acting on $\RG$ which uses the construction of
Heisenberg algebra in \cite{Wa1}. The representation theory of
symmetric groups was earlier studied from the viewpoint of vertex
operators by I.~Frenkel and Jing (cf. \cite{Jing}), and the
natural question of understanding further connections between
vertex algebras and symmetric groups has been posted since then.
Although there has been much nontrivial evidence indicating deep
connections among $\Hx$, $\RG$, and vertex algebras, much remains
to be done to go beyond the vertex representation level in both
pictures\footnote{There has appeared further interesting
development in the symmetric group case by Lascoux and Thibon
\cite{L-T}.}.

The wreath products are relatively simple objects, while the
geometry of Hilbert schemes is very rich. The strikingly parallel
algebraic structures appearing in both $\Hx$ and $\RG$ have been
and will still be a mutual stimulation to study these subjects
together. For instance, the group-theoretic construction of the
Virasoro algebra in \cite{FW} was motivated in part by Lehn's
geometric construction in \cite{Lehn}. The work \cite{FW} in
addition suggested the existence of a surprising relation between
the cup product on the cohomology ring of Hilbert schemes and the
convolution product on the representation ring of wreath products.
A precise connection between the cup product on $(\C^2)^{[n]}$ and
a `filtered' convolution product on $S_n$ has been subsequently
established by Lehn-Sorger \cite{LS} and Vasserot \cite{Vas}. Very
recently Ruan (cf. \cite{Ruan}) observed that such a product
structure on $R(S_n)$ coincides with the orbifold cup product
(introduced in \cite{CR}) for the symmetric product $\C^{2n}/S_n$.
This raises the interesting question how to relate the cohomology
ring structure of Hilbert schemes $\Xn$ to the orbifold cohomology
ring structure of the symmetric products $X^n/S_n$ or more
generally of the wreath product orbifolds $Y^n/\Gn$ \footnote{In
\cite{LS2}, Lehn and Sorger have constructed the cohomology ring
on $\Xn$ in terms of symmetric group $S_n$ and the cohomology ring
$H^*(X)$, when the surface $X$ has numerically trivial canonical
class. The connection between Lehn-Sorger's construction and
Chen-Ruan's orbifold cup product is further clarified in \cite{FG}.}.

In the case when $\G$ is a finite subgroup of $SL_2(\C)$, $Y$ is
$\C^2$, and $X$ is the minimal resolution of singularities $\ale$,
the connection between Hilbert scheme $\alen$ and the wreath
products $\Gn$ is very intriguing, cf. \cite{Wa2}. The natural
morphism $\tau_n: \alen \rightarrow \C^{2n}/\Gn$ is a semismall
crepant resolution of singularities. In \cite{Wa2}, we singled out
a distinguished subvariety $\hquiver$ of the Hilbert scheme
$(\C^2)^{[nN]}$ which is a crepant resolution of $\C^{2n}/\Gn$,
where $N$ is the order of $\G$. It can be shown (cf. the Appendix)
that $\hquiver$ affords a description of quiver varieties in the
sense of Nakajima \cite{Na1, Na4} and it is diffeomorphic to
$\alen$. We established a morphism from the so-called
$\Gn$-Hilbert scheme $\hilqgn$ in $\C^{2n}$ to the $\hquiver$, and
conjectured it is an isomorphism \cite{Wa2} (we refer the reader
to \cite{Nra, Rei1, BKR} for more study of the $G$-Hilbert schemes
in different situations). This conjecture for $n =1$ is a theorem
of Ito-Nakamura \cite{INr1} and Ginzburg-Kapranov (unpublished)
which provides a uniform construction of the minimal resolutions
of simple singularities in terms of Hilbert schemes of points on
$\C^2$. On the other hand, when $\G$ is trivial it is a remarkable
theorem of Haiman \cite{Hai2} which is used to settle the $n!$
conjecture of Garsia-Haiman and the Macdonald polynomial
positivity conjecture \cite{Mac, Hai1}. The combinatorial
implications in our general case (i.e. $\G$ nontrivial including
$\G$ cyclic) along this line will be pursued elsewhere.

The layout of this paper is as follows. In Sect.~\ref{sec_hilb},
we review the generating function for the Betti numbers of Hilbert
schemes, the constructions of Heisenberg algebra and Virasoro
algebra in the framework of Hilbert schemes. We also review
results on the cohomology ring structure of Hilbert schemes. In
Sect.~\ref{sec_wreath}, we explain the constructions of Heisenberg
algebra and Virasoro algebra in the framework of wreath products.
We indicate how to develop further to obtain a group-theoretic
realization of the vertex representations of affine and toroidal
Lie algebras etc. In Sect.~\ref{sec_interplay}, we exhibit the
interplay between the Hilbert schemes and the wreath products. We
discuss the relations between $\alen$ and the $\Gn$-Hilbert scheme
$\hilqgn$. We also review the precise correspondence between the
cup product of $H^*((\C^2)^{[n]})$ and the `filtered' convolution
product of $S_n$. In the Appendix, we present the description of
$\hquiver$ as a quiver variety. In the end, we provide a
dictionary table to make a comparison between Hilbert schemes and
wreath products. We expect such a table will be greatly extended
in the future.

{\bf Acknowledgment}. It is a pleasure to thank Igor Frenkel,
Naihuan Jing, Wei-Ping Li, Zhenbo Qin and Jian Zhou for
stimulating discussions and collaborations on various projects
presented here.

{\bf Note added}. There has been much further development in closely
related topics recently.
We refer the reader to \cite{QW} and the references therein for detail.
\section{Algebraic structures behind Hilbert schemes} \label{sec_hilb}
\subsection{Hilbert schemes $\Xn$}
  Let $X$ be a smooth projective surface over $\C$, and $\Xn$
be the Hilbert scheme of points in $X$. An element in the Hilbert
scheme $\Xn$ is represented by a length-$n$ $0$-dimensional closed
subscheme $\xi$ of $X$. For $\xi \in \Xn$, let $I_{\xi}$ be the
corresponding sheaf of ideals. For a point $x \in X$, let $\xi_x$
be the component of $\xi$ supported at $x$ and $I_{\xi, x} \subset
\mathcal O_{X, x}$ be the stalk of $I_{\xi}$ at $x$. A theorem of
Fogarty \cite{Fog} says that $\Xn$ is smooth. In $\Xn\times X$, we
have the universal codimension-$2$ subscheme:

$$ \mathcal Z_n=\{(\xi, x) \subset \Xn\times X \, | \, x\in
 {\rm Supp}{(\xi)}\}\subset \Xn\times X.
$$

By sending an element in $\Xn$ to its support, we obtain the
Hilbert-Chow morphism $\pi_n: \Xn \rightarrow X^n/S_n$, which is a
resolution of singularities.

For example, when $n=2$, an element of $X^{[2]}$ may be a pair of
two distinct points in $X$, or a pair consisting one point in $X$
together with a tangent direction at this point. The situation
becomes much more complicated for $n>2$.
\subsection{The Betti numbers of Hilbert schemes}
Ellingsrud and Str\o mme \cite{ES1} first calculated the Betti
numbers for the Hilbert scheme $\Xn$ when $X$ is the projective
plane, the affine plane, or a rational ruled surface. They used
the toric action on these surfaces in an essential way, and thus
their method cannot be extended to more general surfaces. Using
the Weil conjecture proved by Deligne, G\"ottsche \cite{Got1,
Got2} calculated the Betti numbers for $\Xn$ associated to an
arbitrary (quasi-)projective surface $X$. Moreover, he presented
the solution in a beautiful generating function which indicates
one should study the cohomology groups of $\Xn$ for all $n\ge 0$
together.

Let us denote by $H^*(-)$ the cohomology group with complex
coefficient, denote by $b_i(-)$ the $i$-th Betti number, and
define the poincar\'e polynomial as
\[
 P_t(X) =\sum_{i \geq 0} t^i b_i(X).
\]

\begin{theorem} \cite{Got1, Got2}
Let $X$ be a smooth quasi-projective surface. The generating
function in a variable $q$ of the Poincar\'e polynomials of the
Hilbert scheme $\Xn$ is given by

\begin{eqnarray*}
 \sum_{n =0}^\infty P_t (\Xn) q^n
 = \prod_{m =1}^\infty \frac
 {(1+t^{2m-1}q^m)^{b_1(X)}(1+t^{2m+1}q^m)^{b_3(X)}}
 {(1-t^{2m-2}q^m)^{b_0(X)}(1-t^{2m}q^m)^{b_2(X)}(1-t^{2m+2}q^m)^{b_4(X)}}
\end{eqnarray*}
\end{theorem}
\
In particular, the above formula implies the following generating
function for the dimension of the total cohomology group $\Hn$:

\begin{eqnarray}   \label{eq_got}
 \sum_{n =0}^\infty \dim \Hn q^n
 = \prod_{m =1}^\infty \frac
 {(1+q^m)^{h^{odd}(X)}}
 {(1-q^m)^{h^{ev}(X)}}
\end{eqnarray}
where $h^{odd}(X) =b_1(X) +b_3(X)$, and $h^{ev}(X) = b_0(X)
+b_2(X) +b_4(X)$. We will denote
\[
\Hx =\bigoplus_{n \ge 0} \Hn.
\]

As remarked by de Cataldo and Migliorini \cite{dCM}, the formulas
above remain valid for any complex surface if one replaces the
Hilbert scheme by a suitable notion called the Douady space.
\subsection{Heisenberg algebra and Hilbert schemes}

Vafa and Witten \cite{VW} observed that the formula (\ref{eq_got})
coincides with the character formula of a Fock space of a
Heisenberg (super)algebra (i.e free bosons/fermions) associated to
the lattice $H^*(X, \Z)/tor$. This motivated them to conjecture
that there indeed exists such a Heisenberg algebra acting on the
direct sum $\Hx = \sum_{n \ge 0} \Hn$ of the cohomology groups of
the Hilbert schemes $\Xn$. For the sake of simplicity, we refrain
ourselves to the case when $X$ is a projective surface when we
recall the construction of Nakajima \cite{Na2} (also compare
\cite{Gro}) below.

For $n \ge 0$ and $\ell \ge 0$, we define $Q^{[n+\ell,n]} \subset
X^{[n+\ell]} \times X \times \Xn$ to be the closed subset $$\{
(\xi, x, \eta) \in X^{[n+\ell]} \times X \times \Xn \, | \, \xi
\supset \eta \text{ and } {\rm Supp} (I_\eta/I_\xi) = \{ x \}
\}.$$ In particular, $Q^{[n,n]} =\emptyset.$

For $n \in \Z$, we define a linear map $\mathfrak q_n: H^*(X)
\mapsto {\rm End} (\Hx)$ as follows. When $n \ge 0$, the linear
operator $\mathfrak q_n(\alpha) \in {\rm End}(\Hx)$ with $\alpha
\in H^*(X)$ is defined by

$$
\mathfrak q_n(\alpha)(a) = \w p_{1*}([Q^{[m+n,m]}] \cdot \w
 \rho^*\alpha \cdot \w p_2^*a).
 $$
for all $a \in H^*(\Xm)$, where $\w p_1, \w \rho, \w p_2$ are the
projections of $\Xmn \times X \times \Xm$ to $\Xmn, X, \Xm$
respectively.

The space $\Hx = \oplus_{n,k} H^k(X^{[n]})$ actually carries a
bi-degree provided by the conformal weight $n$ and the cohomology
degree $k$ respectively. A non-degenerate super-symmetric bilinear
form $(, )$ on $\Hx$ is induced from the standard one on
$H^*(\Xn)$. For a homogeneous linear operator $\mathfrak f \in
\End(\mathbb H)$ of bi-degree $(\ell, m)$, we can define its {\it
adjoint} $\mathfrak f^\dagger \in \End(\mathbb H)$ by $$(\mathfrak
f(\alpha), \beta) = (-1)^{m |\alpha|} \cdot (\alpha, \mathfrak
f^\dagger(\beta))$$ where $|\alpha| = s$ for $\alpha \in H^s(X)$.

When $n < 0$, define the operator $\mathfrak q_n(\alpha) \in {\rm
End}(\Hx)$ with $\alpha \in H^*(X)$ by

$$\mathfrak q_n(\alpha) = (-1)^n \cdot \mathfrak
q_{-n}(\alpha)^\dagger.$$
The operator $\mathfrak q_n (\alpha)$ can be alternatively defined
by switching the role of $\w p_1 $ and $\w p_2$ in the definition
of $\mathfrak q_{-n}(\alpha)$.

\begin{theorem} \cite{Na3}
Let $X$ be a smooth projective surface. The operators $\mathfrak
q_n(\alpha)$ ($n \in \Z$, $\alpha \in H^*(X)$) acting on $\Hx$
satisfy the Heisenberg algebra commutation relations:

$$ [\mathfrak q_n(\alpha), \mathfrak q_m(\beta)] = n \cdot
\delta_{n+m}
 \cdot \int_X(\alpha \beta) \cdot \Id_{\mathbb H}.
$$
Moreover, the space $\Hx = \sum_{n \ge 0} \Hn$ is an irreducible
representation of this Heisenberg algebra.
\end{theorem}

The commutator above is understood as a supercommutator when both
$\alpha$ and $\beta$ are cohomology classes of odd degrees. Indeed
the theorem is valid for an arbitrary quasi-projective surface $X$
with some appropriate modification. We refer to \cite{Na3} for a
proof.
\subsection{Virasoro algebra and Hilbert schemes}

Define the linear operator $\mathfrak d \in \End(\Hx)$ to be the
cup product with $ c_1(p_{1*}\mathcal O_{\mathcal Z_n}) =
-[\partial \Xn]/2$ in $\mathbb H_n$ for each $n$. Here $p_1$ is
the projection of $\Xn \times X$ to $\Xn$, $\partial \Xn$ is the
boundary of $\Xn$ consisting of all $\xi \in \Xn$ with
$|\Supp(\xi)| < n$, and $c_1(p_{1*}\mathcal O_{\mathcal Z_n})$
denotes the first Chern class of the rank-$n$ bundle
$p_{1*}\mathcal O_{\mathcal Z_n}$. Below $K_X$ and $c_2(X)$ stand
for the canonical divisor and the second Chern class of $X$
respectively.

\begin{theorem} \cite{Lehn}
Let $X$ be a smooth projective surface. The commutators of
$\mathfrak d$ and the Heisenberg generators $\mathfrak
q_n(\alpha)$ provide us the Virasoro generators (acting on $\Hx$):

$$[\mathfrak d, \mathfrak q_n(\alpha)] = n \cdot \mathfrak
L_n(\alpha) + \frac{n(|n|-1)}{2} \mathfrak q_n(K_X \alpha),
 $$
where the operators $\mathfrak L_n(\alpha)$ satisfy the Virasoro
algebra commutation relations:
$$ [\mathfrak L_n(\alpha), \mathfrak L_m(\beta)] = (n-m) \cdot
\mathfrak L_{n+m}(\alpha \beta) - \frac{n^3-n}{12} \cdot
\delta_{n+m} \cdot \int_X(c_2(X) \alpha \beta) \cdot \Id_{\Hx}.
 $$
\end{theorem}

Let us recall some elementary construction in vertex algebras, cf.
\cite{Bor1, FLM, Kac2}. Let $a(z) = \sum_{n \in \Z} a_{(n)}
z^{n-\Delta}$ be a vertex operator of conformal weight $\Delta$,
that is, a generating function in a formal variable $z$ where
$a_{(n)}$ is an operator acting on $\Hx$ such that $a_{(n)}
(H^*(\Xm)) \subset \mathbb H^*(\Xmn)$. Put $$
 a_+(z)= \sum_{n>0} a_{(n)} z^{n-\Delta}
 \qquad \text{and} \qquad a_-(z) = \sum_{n
 \le 0} a_{(n)} z^{n-\Delta}
$$
 (note our unusual sign convention here on vertex operators). If
$b(z)$ is another vertex operator, we define a new vertex
operator, which is called {\it the normally ordered product} of
$a(z)$ and $b(z)$, to be:
$$
 :a(z)b(z): = a_+(z) b(z) + (-1)^{ab} b(z) a_-(z)
$$
 where $(-1)^{ab}$ is $-1$ if both $a(z)$ and $b(z)$ are odd and
$1$ otherwise. Inductively we can define the normally ordered
product of $k$ vertex operators from right to left by
$$
 :a_1(z)a_2(z) \cdots a_k(z): \,
 = \, : a_1(z) (:a_2(z) \cdots a_k(z):):.
$$

For $\alpha \in H^*(X)$, we define a vertex operator $\alpha(z)$
by putting
$$
 \alpha(z) = \sum_{n \in \Z} \mathfrak q_n(\alpha) z^{n-1}.
$$

Let $k \ge 1$, $n \in \Z$, and $\alpha \in H^*(X)$. Let
$\delta_{k*}: H^*(X) \to H^*(X^k) \cong H^*(X)^{\otimes k}$ be the
linear map induced by the diagonal embedding $\delta_k: X \to
X^k$, and set
$$
 \delta_{k*}\alpha = \sum_j \alpha_{j,1} \otimes \ldots \otimes
 \alpha_{j, k}.
$$
We define the operator $W^k_n(\alpha) \in \End(\Hx)$ to be the
coefficient of $z^{n-k}$ in the vertex operator
$$
 \frac1{k!} \cdot (\delta_{k*}\alpha) (z)
 = \frac1{k!} \cdot \sum_j : \alpha_{j, 1}(z)
 \cdots \alpha_{j, k}(z):.
$$

Note that $W_n^1(\alpha)$ coincides with the Heisenberg generator
$\mathfrak q_n(\alpha)$, and $W_n^2(\alpha)$ provides us the
Virasoro generator $\mathfrak L_n(\alpha)$, cf. \cite{Lehn}.

It is proved in \cite{LQW1} that if the canonical class $K_X$ is
numerically trivial then
\begin{eqnarray} \label{eq_bound}
 \mathfrak d
 = -W^3_0(1_X).
\end{eqnarray}
It is very interesting to compare with Eqn.~(\ref{eq_conv}).

Recall that $\mathcal Z_n$ is the universal codimension $2$
subscheme of $\Xn \times X$, and $p_1$ and $p_2$ are the
projections of $\Xn \times X$ to $\Xn$ and $X$ respectively. Given
a line bundle (i.e. a locally free sheaf of rank 1) $L$ in $X$,
one can define
 $$L^{[n]} := p_{1*}((\mathcal O_{\mathcal Z_n})
   \otimes p_2^*L).
 $$
Since $p_1$ is a flat finite morphism of degree $n$, $L^{[n]}$ is a
locally free sheaf on $\Xn$ of rank $n$. Using the operator
$\mathfrak d$ effectively, Lehn proved a beautiful formula for the
total Chern class of $L^{[n]}$. This formula was conjectured
earlier by G\"ottsche based on some closely related formula of
Grojnowski and Nakajima (see \cite{Na3}).

\begin{theorem} \cite{Lehn}
 Let $L$ be a line bundle on $X$. Then
 \begin{eqnarray}
  \sum\limits_{n \ge 0}  c ( L^{[n]}) z^n
  &= & \exp \Biggl( \sum_{ n \ge 1}
      ( -1)^{ n -1} \frac 1n \, \mathfrak q_n(c(L))z^n \Biggr).
      \label{eq_class}
 \end{eqnarray}
\end{theorem}

\subsection{The cohomology ring structure of $\Xn$}

For $\gamma \in H^s(X)$ and $n \ge 0$, let $G_i(\gamma, n)$ be the
homogeneous component in $H^{s+2i}(\Xn)$ of
\[
 G(\gamma, n) = p_{1*}(\ch(\mathcal O_{\mathcal Z_n}) \cdot
 p_2^*{\rm td}(X) \cdot p_2^*\gamma) \in \Hn .
\]
Here and below we omit the Poincar\'e duality used to switch a
homology class to a cohomology class and vice versa. The following
theorem has been established by Li, Qin and the author.

\begin{theorem} \cite{LQW1}
Let $X$ be a smooth projective surface. The cohomology ring of the
Hilbert scheme $\Xn$ is generated by the classes $G_{i}(\gamma,
n)$ where $0 \le i < n$ and $\gamma$ runs over a linear basis of
$H^*(X)$.
\end{theorem}

Let us briefly comment on the proof of the above theorem. We took
the viewpoint effectively that for a fixed $\gamma$ the cup
products with $G_{i}(\gamma, n)$ for all $n$ should be treated as
a single operator $\mathfrak G_i(\gamma)$ acting on $\Hx =\oplus_n
\Hn$. Exploring relations between the operators $\mathfrak
G_i(\gamma)$ and the Heisenberg operators $\mathfrak q_n(\alpha)$,
we found some surprising connections between $\mathfrak
G_i(\gamma)$ and the zero-mode of a vertex operator $W^{i+2}(\g)$.
This allowed us to finish the proof by using an induction.

A set of ring generators for $\Hn$ when $X$ is $\mathbb P^2$
(which is easily shown to be equivalent to the set given in the
above Theorem) or $\C^2$ was first found by Ellingsrud and Str\o
mme \cite{ES2}. Lehn's new proof \cite{Lehn} using the Heisenberg
algebra for this result when $X =\C^2$ has been very inspiring for
our approach in \cite{LQW1}. On the other hand, the approach of
Ellingsrud and Str\o mme has been extended to other rational
surfaces and ruled surfaces in \cite{Bea}, and to $K3$ surfaces by
Markman \cite{Mar}.

The generators $G_i(\gamma, n)$ has the advantage of having
certain nice commutation relations with the Heisenberg generators.
But their geometric meaning is not always very clear. A new set of
ring generators for $\Xn$ associated to an arbitrary projective
surface was found in \cite{LQW2} which admits simple geometric
interpretation.

Let us introduce some notations. Let $|0\rangle$ denote the
element $1$ of $H^0(X^{[0]})=\mathbb C$. For $0 \le i < n$ and
$\gamma \in H^i(X)$, define a cohomology class $B_i(\gamma, n) \in
H^{s+2i}(\Xn)$ by putting

$$B_i(\gamma, n) = \frac1{(n-i-1)!} \cdot \mathfrak
q_{i+1}(\gamma) \mathfrak q_{1}(1_X)^{n-i-1} |0\rangle.
 $$
Note that these are among the simplest cohomology classes in
$H^*(\Xn)$ in geometric terms. Indeed, $B_0(1_X, n) = n \cdot
1_{\Xn}$. If either $i > 0$ or $\gamma \in H^s(X)$ with $s > 0$,
then one can easily show that the Poincar\' e dual of $B_i(\gamma,
n)$ is the homology class represented by the closed subset:

$$\{\, \xi\in \Xn\, |\, \hbox{$\exists \, x\in \Upsilon$ with
$\ell(\xi_x)\ge i+1$}\,\}
 $$
where $\Upsilon$ is a homology cycle of $X$ representing the
Poincar\'e dual of $\gamma$, and $\xi_x$ is the component of $\xi$
such that $\xi_x$ is supported at $x$.

\begin{theorem} \cite{LQW2}
Let $X$ be a smooth projective surface. The cohomology ring of the
Hilbert scheme $\Xn$ is generated by the classes $B_{i}(\gamma,
n)$ where $0 \le i < n$ and $\gamma$ runs over a linear basis of
$H^*(X)$.
\end{theorem}

The proof of this theorem is based on establishing certain
relations between $B_{i}(\gamma, n)$ and $G_{j}(\alpha, n)$
introduced earlier and an induction. We also took the viewpoint
that for a fixed $\gamma$ the cup products with $B_{i}(\gamma, n)$
for all $n$ should be treated as a single operator $\mathfrak
B_i(\gamma)$ acting on $\Hx =\oplus_n \Hn$.
\section{Algebraic structures behind wreath products} \label{sec_wreath}
\subsection{The wreath product $\Gn$}

Let $\G$ be a finite group and denote by $\G_*$ the set of
conjugacy classes of $\G$. Denote by $[g]$ the conjugacy class of
$g \in \G$. Denote by $\zeta_c$ the order of the centralizer of an
element $g$ in the conjugacy class $c$ of $\G$. Let $\G^n = \G
\times \ldots \times \G$ be the direct product of $n$ copies of
$\G$. The symmetric group $S_n$ acts on the product group $\G^n$
by permuting the $n$ factors: $ s (g_1, \ldots, g_n) =
(g_{s^{-1}(1)} , \ldots, g_{s^{-1}(n)} )$. The {\em wreath
product} $\Gn$ is defined to be the semi-direct product $\G^n
\rtimes S_n$ of $\G^n$ and $S_n$, namely the multiplication on
$\Gn$ is given by $(g, s)(h, t) = (g. s(h), st)$, where $g, h \in
\G^n, s, t \in S_n$. Note when $\G$ is the trivial one-element
group the wreath product $\Gn$ reduces to $S_n$, and when $\G$ is
$\Z_2$ the wreath product $\Gn$ is the hyperoctahedral group, the
Weyl group of type $B$.

We denote by $R_\Z (\G)$ the lattice with a basis given by the
irreducible characters of $ \G$. Then $R(\G): =R_\Z
(\G)\otimes_{\Z} \C$ can be regarded as the representation ring of
$\G$ or the space of class functions on $\G$. Denote by
 $$\RG =\bigoplus_{n=0}^\infty R(\Gn). $$

The usual bilinear form on $R(\G )$ is defined by

\begin{eqnarray*}
\langle f, g \rangle \equiv \langle f, g \rangle_{\G} = \frac1{ |
\G |}\sum_{x \in \Gamma}
          f(x) g(x^{ -1}),
\end{eqnarray*}
and the set $\G^*$ of the irreducible characters of $\G$ form an
orthonormal basis for $R(\G)$.This provides us a bilinear form
$\langle -, - \rangle$ on $\RG$.

It is well known that the conjugacy classes of $S_n$ are
parameterized by the partitions of $n$. The conjugacy classes of a
general wreath product can be described as follows, (cf. e.g.
\cite{Mac}). Given $a = (g, s) \in \Gn$ where $g = (g_1, \ldots,
g_n)$, we write $s \in S_n$ as a product of disjoint cycles: if
$z= (i_1, \ldots, i_r)$ is one of them, the {\em cycle-product}
$g_{i_r} g_{i_{r-1}} \ldots g_{i_1} $ of $a$ corresponding to the
cycle $z$ is determined by $g$ and $z$ up to conjugacy. For each
$c \in \G_*$ and each integer $r \geq 1$, let $m_r (c)$ be the
number of $r$-cycles in $s$ whose cycle-product lies in $c$.
Denote by $\rho (c)$ the partition having $m_r (c)$ parts equal to
$r$ ($r \geq 1$) and denote by $\rho = ( \rho (c) )_{c \in \G_*}$
the corresponding partition-valued function on $\G_*$. Note that
$|| \rho || : = \sum_{c \in \G_*} |\rho (c)| = \sum_{c \in \G_*, r
\geq 1} r m_r (c) = n$, where $| \rho(c)|$ is the size of the
partition $\rho (c)$. Thus we have defined a map from $\Gn$ to
${\mathcal P}_n (\G_*)$, the set of partition-valued function
$\rho =(\rho(c))_{c \in \G_*}$ on $\G_*$ such that $|| \rho ||
=n$. The function $\rho$ or the data $\{m_r(c) \}_{r,c}$ is called
the {\em type} of $a = (g, s) \in \Gn$. Denote ${\mathcal P}
(\G_*) = \sum_{n \geq 0} {\mathcal P}_n (\G_*)$. It can be shown
that two elements in $\Gn$ are conjugate to each other if and only
if they have the same type.
\subsection{Heisenberg algebra and wreath products}
Denote by $c_n (c \in \G_*)$ the conjugacy class in $\Gn$ of
elements $(x, s) \in \Gn$ such that $s$ is an $n$-cycle and the
cycle product of $x$ is $ c$. Denote by $\sigma_n (c )$ the class
function on $\Gn$ which takes value $n \zeta_c$ (i.e. the order of
the centralizer of an element in the class $c_n$) on elements in
the class $c_n$ and $0$ elsewhere. Given $\gamma \in R(\G)$, we
denote by $\sigma_n (\g )$ the class function on $\Gn $ which
takes value $n \g (c) $ on elements in the class $c_n, c \in
\G_*$, and $0$ elsewhere. In the symmetric group (i.e. $\G$
trivial) case, the definition of $\sigma_n(1)$ is much simplified.

We define $ p_{n} (\gamma), n >0$ to be a map from $\RG$ to itself
by the following composition:
\[  R (\Gm) \stackrel{ \sigma_n ( \g ) \otimes }{\longrightarrow}
  R(\Gn) \bigotimes R (\Gm)  \stackrel{{Ind} }{\longrightarrow}
  R ( {\Gamma}_{n +m}).
\]
We define $ p_{-n} (\gamma), n >0$ to be the adjoint of
$p_n(\gamma)$ with respect to the bilinear form $\langle -, -
\rangle$ on $\RG$, or equivalently a map from $\RG$ to itself as
the composition
\[
  R (\Gm)  \stackrel{ Res }{\longrightarrow}
   R(\Gn)\bigotimes R ( {\G }_{m -n})
 \stackrel{ \langle \sigma_n ( \g), \cdot \rangle}{\longrightarrow}
 R ( {\G }_{m -n}).
\]

The following theorem is proved by using Mackey's theorem on the
restriction to a subgroup of an induced character. Indeed it is a
special case of a theorem established in \cite{Wa1} in a more
general (i.e. equivariant $K$ groups) setup.

\begin{theorem}  \cite{Wa1}
The operators $p_n(\gamma)$, $n \in \Z, \gamma \in R(\G)$, satisfy
the Heisenberg commutation relations:

$$[p_m( \gamma), p_n(\gamma ')]
 = -m \delta_{m, -n} \langle \gamma, \gamma' \rangle C,
 \quad \g, \g ' \in R(\G).
 $$
Moreover, $\RG$ is an irreducible representation of the Heisenberg
algebra.
\end{theorem}

It follows that $\RG$ is isomorphic to the symmetric algebra
generated by $p_m(\g)$ where $ \g \in \G^*, n>0$. The construction
of Heisenberg algebra above is intimately related to the Hopf
algebra structure on $\RG$ given by Zelevinsky \cite{Zel}. One may
choose a different natural basis of the Heisenberg algebra
parameterized by $n \in \Z, c \in \G_*$ (cf. \cite{FJW1}):
\[
 p_n( c) = \sum_{ \g\in \G^*} \gamma(c^{-1}) p_n( \g ),
\]
where $c^{-1}$ denotes the conjugacy class $\{g \in \G| \g^{-1}\in
c \}$.

Given a representation $V$ of $\G$ with character $\g$, the
natural actions of $\G^n$ and $S_n$ on $n$-th outer tensor product
$V^{ \otimes n} $ of $V$ can be combined into an action of the
wreath product $\Gn$. We denote the $\Gn$-character of $V^{
\otimes n} $ by $\eta_n (\g)$. Denote by $\varepsilon_n$ the
(1-dimensional) sign representation of $\Gn$ on which $\G^n$ acts
trivially while $S_n$ acts as sign representation. We denote by
$\varepsilon_n ( \g ) \in R(\Gn)$ the character of the tensor
product of $\varepsilon_n$ and $V^{\otimes n}$.

The following proposition can be found in \cite{Wa1} in a more
general (i.e. equivariant $K$-group) setup. These equalities here
are well known for $\G \in \G^*$, cf. \cite{Mac, FJW1}. We remark
that this type of formulas together with a parallel formula
(\ref{eq_class}) often appear as half of a vertex operator in
vertex algebra literature.

\begin{proposition}  \label{prop_exp}
   For any $\g \in R(\G)$, we have
\begin{eqnarray}
 \sum\limits_{n \ge 0}  \ch ( \eta_n( \g ) ) z^n
  &= & \exp \Biggl( \sum_{ n \ge 1}
      \frac 1n \, a_{-n}(\g )z^n \Biggr),
      \nonumber 
      \\
 \sum\limits_{n \ge 0}  \ch ( \varepsilon_n( \g )  ) z^n
  &= & \exp \Biggl( \sum_{ n \ge 1}
      ( -1)^{ n -1} \frac 1n \, a_{-n}(\g )z^n \Biggr). \label{eq_sign}
\end{eqnarray}
\end{proposition}
\subsection{McKay correspondence and wreath products}

In the case when $\G$ is a finite subgroup of $SL_2(\C)$, it is
pointed out in \cite{Wa1} that one can indeed develop the above
construction further to realize the basic vertex representations
of affine Lie algebras of ADE types in a group theoretic manner.
This has subsequently been fully carried out in collaboration with
I.~Frenkel and Jing \cite{FJW1}.

The classification of finite subgroups of $SL_2 (\C)$ is of course
well known. First it is clear that finite subgroups of $SL_2 (\C)$
are indeed finite subgroups of $SU_2$. Note that $SU_2$ is a
double cover of $SO_3$, so the classification is essentially
reduced to the classification of finite subgroups of $SO_3$, which
consists of the symmetric groups of regular polyhedra (i.e.
tetrahedral, octahedral and icosahedral groups) and the symmetric
groups of certain two dimensional `degenerate regular polyhedra'
(i.e. cyclic and dihedral groups). This results into a complete
list of finite subgroups of $SL_2(\C)$: the cyclic, binary
dihedral, tetrahedral, octahedral and icosahedral groups. (Binary
here means a double cover).

Denote by $\g_0$ the trivial character of $\G$ and $Q$ the
two-dimensional defining representation of $\G$ in $\C^2$. Let
$\wt = 2\g_0 -Q$. Introduce the {\em weighted} bilinear form (cf.
\cite{FJW1}) on $R(\G)$ by
\[
 \langle f, g \rangle_{\wt } = \langle \wt \otimes f ,  g \rangle_{\G },
   \quad f, g \in R( \G).
\]
Then McKay's observation \cite{McK} simply says the lattice
$(R_\Z(\G), \langle -, - \rangle_{\wt })$ is the affine root
lattice of ADE types, and this gives rise to a bijection between
the set of finite subgroups of $SL_2(\C)$ and the set of affine
root lattice (or equivalently affine Dynkin diagrams) of ADE
types. For example, cyclic corresponds to type $A$.

Recall $\eta_n ( \wt )$ is defined to be a virtual character in
$R(\Gn)$. We define a weighted bilinear form on $R(\Gn)$ by
\[
 \langle f, g \rangle_{\wt } =\langle\eta_n(\wt)\otimes f ,g \rangle_{\Gn},
   \quad f, g \in R( \Gn),
\]
and thus a bilinear form $\langle -, - \rangle_{\wt }$ on $\RG
=\oplus_n R(\Gn)$.

Let us digress for a moment. Given an integral lattice $L$, we can
construct a Heisenberg algebra $H_L$ which is the central
extension of the loop algebra $\mathfrak h \otimes \C [t,t^{-1}]$,
where $\mathfrak h$ denotes $L\otimes_\Z \C$. One can endow a
canonical structure of an irreducible representation of $H_L$ on
the symmetric algebra $\mbox{Sym}(\mathfrak h \otimes t^{-1}
\C[t^{-1}])$. Following the usual consideration for a lattice
vertex algebra \cite{Bor1, FLM, Kac2}, we introduce the space
$V_L =\mbox{Sym}(\mathfrak h \otimes t^{-1} \C[t^{-1}]) \bigotimes
\C[L],$
where $\C[L]$ denotes the group algebra of the lattice (i.e. free
abelian group) $L$. A celebrated theorem of Frenkel-Kac \cite{FK}
and Segal \cite{Seg1} says that if $L$ is a simple root lattice of
ADE type then $V_L$ provides an realization of the basic
representation of the affine Kac-Moody algebra (which is the
central extension of a loop algebra associated to the lattice
$L$), cf. \cite{Kac1}. The construction is based on the action of
the Heisenberg algebra $H_L$ and vertex operators. This
construction has been extended in \cite{MRY} to realize the vertex
representation of a toroidal Lie algebra (which is a central
extension of a double loop algebra associated to $L$) in the case
when $L$ is an affine root lattice.

We now introduce the following space
\[
 V_\G = \RG \bigotimes \C[R_\Z(\G)],
\]
where $\C[R_\Z(\G)]$ is the group algebra of the lattice
$(R_\Z(\G), \langle -, - \rangle_{\wt })$.

We formulate informally a main theorem in \cite{FJW1} below, and
refer the reader to the original paper for more detail. Our
constructions can be regarded as a new form of McKay
correspondence: starting from a finite subgroup $\G$ of
$SL_2(\C)$, we realize the vertex representation of an affine and
toroidal Lie algebra whose Dynkin diagram corresponds to $\G$ in
the classical way.

\begin{theorem}
For a finite subgroup $\G$ of $SL_2(\C)$, the space $V_\G$
provides a group-theoretic realization of the vertex
representation of a toroidal Lie algebra associated to the affine
root lattice $(R_\Z(\G), \langle -, - \rangle_{\wt })$. The
generators for the toroidal Lie algebra are given in terms of
vertex operators and afford a natural group-theoretic
interpretation in terms of certain induction/restriction functors.

One can identify a distinguished subspace $\overline{V}_\G$ of
$V_\G$ (which is, roughly speaking, the subspace to which the
trivial character $\g_0$ of $\G$ makes no contribution) in which
we realize in a group-theoretic manner the basic vertex
representation of an affine Lie algebra.
\end{theorem}

Note that the action of $\Gn$ on $\C^{2n}$ induces a natural
degree-preserving $\Gn$-action on the exterior algebra $ \Lambda^*
(\C^{2n}) = \oplus_{i=0}^{2n} \Lambda^i (\C^{2n}). $ It we
understand $\C^{2n}$ as the $\Gn$-equivariant vector bundle over a
point, then $\Lambda^* (\C^{2n})$ is exactly the associated
Koszul-Thom complex. Note that the virtual $\G$-character $\wt$
can be rewritten as $\sum_{i=0}^{2} (-1)^i \Lambda^i (\C^{2})$.

\begin{proposition} \cite{Wa2}
The virtual $\Gn$-character $\sum_{i=0}^{2n} (-1)^i \Lambda^i
(\C^{2n})$ coincides with the virtual $\Gn$-character
$\eta_n(\wt)$, where the $\G$-virtual character $\wt$ is given by
$\sum_{i=0}^{2} (-1)^i \Lambda^i (\C^{2})$.
\end{proposition}

\begin{remark} \rm
We denote by $D_{\Gn} (\C^{2n})$ the bounded derived category of
$\Gn$-equivariant coherent sheaves on $\C^{2n}$. We denote by
$D^0_{\Gn} (\C^{2n})$ the full subcategory of $D_{\Gn} (\C^{2n})$
consisting of objects whose cohomology sheaves are concentrated on
the origin of $\C^{2n}$, and by $K^0_{\Gn} (\C^{2n})$ the
corresponding $K$-group. A natural bilinear form on $K^0_{\Gn}
(\C^{2n})$ uses the $\Gn$-virtual character $\sum_{i=0}^{2n}
(-1)^i \Lambda^i (\C^{2n})$. It follows from the above proposition
that there exists a canonical isometric isomorphism between
$K^0_{\Gn} (\C^{2n})$ endowed with this bilinear form and
$R_{\Z}(\Gn)$ endowed with the weighted bilinear form, cf.
\cite{Wa2}.
\end{remark}

The realization of the homogeneous vertex representation of an
affine algebra on $\overline{V}_\G$ above can be regarded as a
counterpart of a realization in terms of the middle dimensional
cohomology group of the Hilbert schemes $\Xn$ when $X$ is the
minimal resolution $\ale$ of the simple singularity $\C^2/\G$, cf.
\cite{Na3}. However the toroidal Lie algebras do not appear in
such a cohomology group setup since one cannot put the cohomology
groups of different degrees of $\Xn$ on the same footing. We
believe if one can construct the Heisenberg algebra on the direct
sum of the {\rm $K$-groups} of $\alen$ for all $n \ge 0$, then one
would be able to obtain toroidal algebras as well. It would be
very interesting to establish a {\rm direct} isomorphism between
these two constructions of toroidal Lie algebras.

By adding $\C^*$ to the picture appropriately and noting the
representation ring of the algebraic representations of $\C^*$ can
be identified with the polynomial ring $\C[q, q^{-1}]$, we
\cite{FJW2} obtained a group-theoretic realization of the
Frenkel-Jing construction of quantum affine algebras and a
construction of the quantum toroidal algebras by Ginzburg,
Kapranov and Vasserot. We remark that the fact that the
two-dimensional defining representation $Q$ of $\G$ is reducible
for $\G$ cyclic is reflected by the fact that the toroidal algebra
of type $A$ affords a {\em two}-parameter deformation. It would be
interesting to obtain a realization of vertex representations of
quantum toroidal algebras using the equivariant $K$-groups
$K_{\C^*}(\alen)$.

Using spin/projective representations of wreath products, we have
realized in \cite{FJW3, JW} the vertex representations of certain
twisted affine and toroidal algebras.
\subsection{Virasoro algebra and wreath products}

Given $f, g \in \C [\G]$, the  {\em convolution} product on $
\C[\G]$ is defined by
\[
  (f * g) (x) =
    \sum_{y \in \G} f(x y^{-1}) g(y),
  \quad f, g \in \C [\G], x \in \G.
\]
In particular if $f, g \in R(\G)$, then so is $f * g$. It is well
known that
\begin{eqnarray}  \label{eq_idem}
   \g ' * \g = \frac{|\G|}{d_{\g}} \delta_{\g ,\g '} \g,
    \quad \g', \g \in \G^*,
\end{eqnarray}
where $d_{\g}$ is the {\em degree} of the irreducible character
$\g$.

Denote by $K_c$ the sum of all elements in a conjugacy class $c$.
By abuse of notation, we also regard $K_c$ the class function on
$\G$ which takes value $1$ on elements in the conjugacy class $c$
and $0$ elsewhere. It is clear that $K_c, c \in \G_*,$ form a
basis of $R(\G)$. The elements $K_c, c \in \G_*$ actually form a
linear basis of the center in the group algebra $\C[\G]$.

Given $c \in \G_*$, we denote by $K_i (c,n)$ the conjugacy class
corresponding to the partition-valued function which maps the
class $c$ to the one-part partition $(i+1)$, $c^0$ to the
partition $(1^{n-i-1})$ and other classes in $\G_*$ to $0$. Let us
define an operator $\mathfrak \Delta_i(K_c)$ on $\RG$ to be the
convolution product with $K_i(c,n)$ in $R(\Gn)$ for each
$n$.
The $K_i(c,n)$
is the counterpart of the $B_i(\g,n) \in \Hn$
studied in Section~\ref{sec_hilb}. In the symmetric group case,
I.~Frenkel and I had much discussion on $K_i(c,n)$ in early 1999.

\begin{theorem} \cite{FW}
The commutation relation between the operator $\mathfrak
\Delta_1(c)$ associated to a conjugacy class $c \in \G_*$ and the
Heisenberg algebra generator $p_{n} (\g)$ is given by
 \[
 [\mathfrak \Delta_1(K_c), p_{n} (\g)]
  = \frac{n|\G |^2 \g (c^{-1})}{\zeta_c d_{\g}^2} L_{n}(\g),
\]
where the operators $L_n (\g)$ acting on the space $\RG$ satisfy
the Virasoro commutation relations:
 \[
    [ L_n(\g), L_m(\g ') ]
    = (n-m) \delta_{\g, \g'}L_{n +m}(\g) -
     \frac{n^3 -n}{12} \delta_{\g, \g'} \delta_{n, -m},
 \]
where $\g, \g '$ lie in $\G^*$.
\end{theorem}

We can extend the operator $\mathfrak \Delta_1(K_c)$ linearly to any $a \in
R(\G)$. We remark that the following {\em transfer property}
holds: $[\mathfrak \Delta_1(a), p_{n} (b)]$, where $a, b \in
R(\G)$, does not depend on $a$ and $b$ individually but only on
the product $a*b$. Similar transfer properties have been observed
for the various geometric operators introduced in the framework of
Hilbert schemes \cite{Lehn, LQW1, LQW2}, cf. Sect.~\ref{sec_hilb}.

In the symmetric group case, the operator $\mathfrak \Delta_1
=\mathfrak \Delta_1(c^0)$ has been studied by Goulden \cite{Gou}
for a totally different purpose and rediscovered in \cite{FW} in
the search for the Virasoro algebra. In this case we have
\begin{eqnarray}
 \mathfrak \Delta_1 &=& \frac16 Res_{z=0} : a (z)^3: z^2
 \nonumber \\
 &=& \frac12 \sum_{n, m> 0}
  ( p_n p_m p_{ -n -m} + p_{n+m} p_{-n}p_{-m}), \label{eq_conv}
\end{eqnarray}
where $: a (z)^3:$ denotes the normally ordered product of $a(z) =
\sum_{n \in  \Z}p_n z^{-n-1}$ with itself for three times. In
general, $\mathfrak \Delta_1(c)$ has a similar expression as an operator
acting on $\RG$, cf. \cite{FW} for more detail.
\section{Interplay between Hilbert schemes and wreath products}
\label{sec_interplay}
\subsection{Desingularization of the wreath product orbifolds}
 Now we assume that $Y$ is a quasi-projective surface acted upon by
a finite group $\G$, and assume that a resolution of singularities
$\tau: X \rightarrow Y/\G$ is given. The action of the product
group $\G^n$ and the symmetric group $S_n$ on the direct product
$Y^n$ induces a natural action of the wreath product $\Gn$ on
$Y^n$: for $ a = ( (g_1, \ldots, g_n), s)$, we let
\begin{eqnarray*}
 a . (x_1, \ldots, x_n)
  = (g_1 x_{s^{-1} (1)}, \ldots, g_n x_{s^{-1} (n)})
\end{eqnarray*}
where $x_1, \ldots, x_n \in Y$. It is observed in \cite{Wa1} that
the following commutative diagram

\begin{eqnarray} \label{diag_resol}
\CD \Xn @>{\pi_n}>> X^n  /S_n \\ @VV{\tau_n}V @VV{\tau_{(n)}}V
\\ Y^n /\Gn @<{\cong}<< (Y/\G)^n /S_n
\endCD
\end{eqnarray}
defines a resolution of singularities $\tau_n: \Xn \rightarrow
Y^n/ \Gn$. Here we used the observation that the wreath product
orbifold is exactly the symmetric product of an orbifold:
$(Y/\G)^n /S_n = (Y^n/\G^n) /S_n =Y^n /(\G^n \rtimes S_n) =Y^n
/\Gn$.

\begin{remark} \rm
A direct sum $\mathcal F_\G(Y)$ of equivariant $K$-groups
$K_{\Gn}(Y^n) \otimes \C$ for $n \ge 0$ is shown \cite{Wa1} to
have several nice algebraic structures. In particular, a
Heisenberg algebra associated to the lattice $K_\G(Y)$ modulo
torsion, constructed in terms additive $K$-group maps, acts on
$\mathcal F_\G (Y)$ irreducibly. This specializes to the
construction of Segal \cite{Seg1} for $\G$ trivial, and
specializes to the construction on $\RG =\oplus_n R(\Gn)$ in
Sect.~\ref{sec_wreath} when $Y$ is a point.
\end{remark}
\subsection{The resolution $\Xn \rightarrow Y^n/\Gn$ is good}\label{subsec_goodres}

A guidline here \cite{Wa1} is that whenever the resolution $\tau:
X \rightarrow Y/\G$ is `good' \/ in a suitable sense then the
resolution $\tau_n: \Xn \rightarrow Y^n/\Gn$ is also `good';
alternatively, whenever a reasonable statement can be made relating the
resolution $X$ and the orbifold $Y/\G$, then the corresponding
statement is true relating the resolution $\Xn$ and the orbifold
$Y^n/\Gn$.

For example, if $\tau: X \rightarrow Y/\G$ is a crepant,
respectively semismall, resolution of singularities, then we can
show that $\tau_n: \Xn \rightarrow Y^n/\Gn$ is also a crepant,
respectively semismall, resolution of singularities.

Let us recall that in the study of orbifold string theory Dixon
{\em et al} \cite{DHVW} introduced the {\em orbifold Euler number}
as

$$\chi (M, G) = \frac1{|G|} \sum_{g_1 g_2 = g_2 g_1} \chi(M^{g_1,
g_2}),$$ where $M$ is a smooth manifold acted by a finite group
$G$, $M^{g_1, g_2}$ denotes the common fixed point set of $g_1$
and $g_2$ and $\chi (\cdot)$ denotes the usual Euler number.

One can also define the orbifold Hodge numbers $h^{s,t}(M, G)$ if
$M$ is a compact complex manifold (cf. Zaslow \cite{Zas}). For
each conjugacy class $c =[g] \in G_*$, let $Y^g_1, \cdots,
Y^g_{N_c}$ be the connected components of the fixed-point set
$Y^g$. On the tangent space to each point in $Y^g_{\alpha}$, $g$
acts as a diagonal matrix $\diag(e^{ 2\pi\sqrt{-1} \theta_1},
\cdots, e^{2\pi\sqrt{-1}\theta_d})$, where $0 \leq \theta_i <1$.
The {\rm shift} number $F^g_{\alpha} =\sum_{j=1}^d \theta_j$ is an
integer if we further assume $g$ acts on the tangent space by a
matrix in $SL(n, \C)$. The {\em orbifold Hodge numbers} of the
orbifold $M/G$ are defined to be
\begin{eqnarray*}
 h^{s,t}(Y, G) = \sum_{c =[g] \in G_*}
\sum_{\alpha_c =1}^{N_c} h^{s- F^c_{\alpha_c},
t-F^c_{\alpha_c}}(Y^g_{\alpha_c}/Z(g)),
\end{eqnarray*}
where $Z(g)$ is the centralizer of $g$ in $G$, and
$h^{*,*}(Y^g_{\alpha_c}/Z(g))$ denotes the dimension of the space
of invariants $H^{*,*}(Y^g_{\alpha_c})^{Z(g)}$. We define the
orbifold virtual Hodge polynomial by $e(M,G;x,y) =\sum_{s,t}
(-1)^{s+t}h^{s,t}(M,G) x^s y^t.$ In particular when $G$ is
trivial, the orbifold Euler number, orbifold Hodge numbers and
orbifold virtual Hodge polynomials reduce to the usual ones.

Now let us get back to our setup when $X,Y$ are complex surfaces,
and $\tau: X \rightarrow Y/\G$ is a resolution of singularities.
The Euler and Hodge numbers of Hilbert schemes were calculated in
\cite{Got1, Got2, GS} (also cf. \cite{Che}). On the other hand,
the orbifold Euler and Hodge numbers for the wreath product
orbifolds were calculated in \cite{Wa1, WaZ}. In terms of
generating functions, we have

\begin{eqnarray*}
 \sum_{n =0}^\infty \chi(\Xn) q^n
 &=& \prod_{m =1}^\infty (1-q^m)^{- \chi(X)}, \\
 \sum_{n =0}^\infty \chi(Y^n, \Gn) q^n
 &=& \prod_{m =1}^\infty (1-q^m)^{- \chi(Y, \G)}, \\
 %
 \sum_{n=1}^{\infty} e(\Xn; x,y) q^n
 &=& \prod_{r=1}^{\infty}\prod_{s,t} (1 - x^sy^t
 q^r(xy)^{(r-1)})^{(-1)^{s+t+1} h^{s, t}(X)}, \\
 \sum_{n=1}^{\infty} e(Y^n, \Gn; x,y) q^n
 &=& \prod_{r=1}^{\infty}\prod_{s,t} (1 - x^sy^t
 q^r(xy)^{(r-1)})^{(-1)^{s+t+1} h^{s, t}(Y, \G)}.
\end{eqnarray*}

As observed in \cite{Wa1, WaZ}, an immediate consequence of the
above formulas is that if the Euler (resp. Hodge) number of $X$ is
equal to the orbifold Euler (resp. Hodge) number of the orbifold
$Y/\G$, then the Euler (resp. Hodge) number of $\Xn$ is equal to
the orbifold Euler (resp. Hodge) number of the orbifold $Y^n/\Gn$.
To our best knowledge, the resolutions $\tau_n:\Xn \rightarrow
Y^n/\Gn$ cover all known higher dimensional examples which match
the Euler and Hodge numbers between the resolution and the
orbifold. It has been conjectured \cite{WaZ} that a similar
statement holds by substituting (orbifold) Hodge numbers with
(orbifold) elliptic genera (compare \cite{DMVV}). Recently Borisov
and Libgober \cite{BL} obtained a mathematical formulation for
orbifold elliptic genera and verified a formula for the elliptic
genera of $Y^n/\Gn$ conjectured in \cite{WaZ}. As we can easily
construct various examples of good resolutions $\tau: X
\rightarrow Y/\G$, our construction provides in a tautological way
numerous higher dimensional examples of good resolutions.

Note that the guiding principle, when applied to the situation
when $\G$ is trivial and $X$ equals $Y$, simply says that the
Hilbert-Chow morphism $\Xn \rightarrow X^n  /S_n$ is a `good'
resolution of singularities in every aspect, cf. \cite{HH, Got3,
Na3, DMVV, Zhou}. Another very interesting example is the $K3$
surfaces which provide good resolutions of the $\Z/2\Z$ quotient
of abelian surfaces. The coincidence between the corresponding
Hodge numbers of Hilbert schemes and wreath product orbifolds when
$X$ is $K3$ has also been obtained independently by Bryan, Donagi,
and Leung \cite{BDL}.

A very important example is as follows. Let $\G$ be a finite
subgroup of $ SL_2 (\C )$. We denote by $\tau : \ale \rightarrow
\C^2 /\G$ the minimal resolution of a simple singularity. By
applying the above general construction to the case $Y = \C^2$ and
$X = \ale$, we obtain the following commutative diagram:

\begin{eqnarray} \label{diag_ale}
\CD \alen @>{\pi_n}>> (\ale)^n  /S_n \\ @VV{\tau_n}V
@VV{\tau_{(n)}}V \\ \C^{2n} /\G_n @<{\cong}<< (\C^2/\G)^n /S_n
\endCD
\end{eqnarray}
This example will be examined from a different viewpoint in the
next subsection.
\subsection{The $\Gn$-Hilbert scheme $\hilqgn$}
 Let $Y$ be a smooth complex algebraic variety and let $\G$ be
a finite subgroup of order $N$ in the group of automorphisms ${\rm
Aut}(Y)$. A regular $\G$-orbit can be viewed as an element in the
Hilbert scheme $Y^{[N]}$ of $N$ points in $Y$. The {\em
$\G$-Hilbert scheme} (or Hilbert scheme of $\G$-regular orbits in
$Y$) $\hilq$ in $Y$ is defined to be the closure of the set of
regular $\G$-orbits in $Y^{[N]}$ (cf. e.g. \cite{Nra, Rei1}).
There is an induced $\G$-action on $Y^{[N]}$, and it is clear that
$\hilq$ is a component of the set $Y^{[N],\G}$ of $\G$-fixed
points in $Y^{[N]}$.

A theorem of Ito-Nakamura \cite{INr1, INr2} (also observed by
Ginzburg and Kapranov) says the $\G$-Hilbert scheme $\ale$
associated to a finite group $\G$ of $SL_2(\C)$ is the minimal
resolution of $\C^2/\G$. So our notation for orbit Hilbert schemes
cause no real confusion with the previous sections where $\ale$ is
used to denote the minimal resolution. We refer the reader to
\cite{INr2} for many other interesting connections with other
branches of mathematics.

The Hilbert-Chow morphism from $Y^{[N]}$ to $Y^{(N)}$ induces a
morphism $\hilq \rightarrow Y/\G$. The notion of orbit Hilbert
schemes is particularly important when the morphism $\hilq
\rightarrow Y/\G$ turns out to be a (crepant) resolution of
singularities (cf. \cite{Rei1, Nra, BKR} and references therein).
But in practice it is a very difficult problem to determine
whether the $\G$-Hilbert scheme is smooth and even more
difficult to have a good description of it.

Below we restrict our study to the setup of the diagram
(\ref{diag_ale}). This material is taken from \cite{Wa2}. We
denote by $\hquiver$ the set of $\G$-invariant ideals $I$ in the
Hilbert scheme $ (\C^2)^{[nN]}$ such that the quotient $\C [x, y]
/ I$ is isomorphic to a direct sum $R^n$ of $n$ copies of the
regular $\G$-representation $R$ as a $\G$-module, where $N$ is the
order of the group $\G$. Since the quotient $\C [x, y]/I$ are
isomorphic as $\G$-modules for all $I$ in a given connected
component of $(\C^2 )^{[nN], \G}$, the variety $\hquiver$ is
non-singular, and is a union of components of $(\C^2 )^{[nN],
\G}$. It turns out that $\hquiver$ is connected as we will explain
later.

Given $J \in \hilqgn$, we regard it as an ideal in $\C [{\bf x},
{\bf y}]$ of length $N^n n!$ (which is the order of $\Gn$), where
${\bf x}$ and ${\bf y}$ stand for the coordinates $x_1,
x_2,\ldots, x_n$ and $y_1, y_2, \ldots, y_n$ in $\C^{2n}$. Then
the quotient $\C [{\bf x}, {\bf y}] / J$ affords the regular
representation of $\Gn$, and its only $\Gn$-invariants are
constants. Thus for a given $f \in \C [x,y]^{\G}$, we have
$\sum_{i=1}^n f(x_i, y_i) = c_f \mbox{ mod }J$ for some constant
$c_f$.

Denote by $\G_{n -1}$ the subgroup of $\Gn$ acting on $\C
[x_2,\ldots, x_n, y_2, \ldots, y_n]$. We can show by induction
that the algebra of invariants $\C [x_2,\ldots, x_n, y_2, \ldots,
y_n]^{\G_{n -1}}$ is generated by the polynomials $\sum_{i=2}^n
f(x_i, y_i)$, where $f$ runs over an arbitrary linear basis
$\mathcal B$ for the space of invariants $\C [x,y]^{\G}$, (cf.
\cite{Wa2}). In the case when $\G$ is trivial, this is a theorem
of Weyl. It follows that the space $\C [{\bf x}, {\bf y}]^{\G_{n
-1}}$ is generated by $x_1, y_1$ and $\sum_{i=2}^n f(x_i, y_i)$,
where $f \in \mathcal B$. The latter is equal to $c_f - f(x_1,
y_1) \mbox{ mod } J$. Thus $(\C [{\bf x}, {\bf y}] / J)^{\G_{n
-1}}$ is generated by $x_1, y_1$ and $c_f - f(x_1, y_1)$, where $f
\in \C [x, y]^{\G}$. It follows that
\begin{eqnarray}  \label{eq_link}
  \C [x_1, y_1 ] / (J \cap \C [x_1, y_1 ])
  \cong ( \C [{\bf x},{\bf y}] / J)^{\G_{n -1}}.
\end{eqnarray}
This space has dimension $nN = | \Gn | / | \G_{n -1} |$ because
$(\C[{\bf x}, {\bf y}] /J)^{\G_{n -1}}$ can be identified with the
space of ${\G_{n -1}}$-invariants in the regular representation of
$\Gn$. The first copy of $\G$ in the Cartesian product $\G^n
\subset \Gn$ commutes with $\G_{n -1}$ above. It follows from
(\ref{eq_link}) that the quotient $\C [x_1, y_1 ] / (J \cap \C
[x_1, y_1 ])$ as a $\G$-module is isomorphic to $R^n$. That is, $J
\cap \C [x_1, y_1 ]$ lies in $\hquiver$.

The map $\varphi$ can be also understood as follows. Let
$\quotuniv$ be the universal family over $\hilqgn$ which is a
subvariety of the Hilbert scheme $(\C^{2n})^{[n!N^n]}$:

\begin{eqnarray*} \label{diag_univ}
 \CD \quotuniv @>>> \C^{2n} \\ @VVV @VVV \\ \hilqgn  @>>>
\C^{2n} /\Gn \endCD
\end{eqnarray*}
It has a natural $\Gn$-action fiberwise such that each fiber
carries the regular representation of $\Gn$. Then $\quotuniv/
\G_{n-1}$ is flat and finite of degree $nN$ over $\hilqgn$, and
thus can be identified with a family of subschemes of $\C^2$ as
above. Then $\varphi$ is the morphism given by the universal
property of the Hilbert scheme $(\C^2 )^{[nN]}$ for the family
$\quotuniv/ \G_{n-1}$.

It is observed \cite{Wa2} that the variety of $\G$-fixed-points
$((\C^2)^{nN}/S_{nN})^{\G}$ can be naturally identified with the
wreath product orbifold $\C^{2n}/\Gn$. Thus the Hilbert-Chow
morphism $(\C^2)^{[nN]} \rightarrow (\C^2)^{nN}/S_{nN}$ induces a
morphism between the varieties of $\G$-fixed-points $\varsigma_n:
\hquiver \rightarrow \C^{2n}/\Gn$. As a $\G$-fixed-point set,
$\hquiver$ has an induced holomorphic symplectic 2-form from the
Hilbert scheme $(\C^2)^{[nN]}$. It is easy to see that over the
set of generic points in $\C^{2n}/\Gn$ which consist of
$\Gn$-regular orbits on $\C^{2n}$, the morphism $\varsigma_n$ is
one-to-one. It follows that $\varsigma_n$ is a crepant resolution
of singularities.

\begin{theorem} \cite{Wa2}
There exists a natural morphism $\varphi : \hilqgn \longrightarrow
\hquiver$ such that the following diagram commutative:
\begin{eqnarray*}
 \CD \hilqgn @>{\varphi}>> \hquiver \\
@V{\sigma_n}VV @V{\varsigma_n}VV
\\  \C^{2n} /\Gn @= \C^{2n} /\Gn \endCD
\end{eqnarray*}
Here $\varsigma_n: \hquiver \rightarrow \C^{2n}/\Gn$ is a crepant
resolution of singularities. In particular $\varphi$ is bijective
over the set of generic points on $\hilqgn$ consisting of the
regular $\Gn$ orbits on $\C^{2n}$.
\end{theorem}

It is of great iterest to clarify the precise relations among
$\hilqgn, \hquiver$ and $\alen$. In an optimistic way, it was
conjectured in \cite{Wa2} that the morphism $\varphi : \hilqgn
\longrightarrow \hquiver$ is an isomorphism over $\C^{2n} /\Gn$.
In the Appendix, we explain that $\hquiver$ admits a quiver
variety description in the sense of Nakajima \cite{Na1, Na3} (cf.
\cite{Wa2}). In particular this implies that $\hquiver$ is
connected. According to Nakajima (unpublished), the Hilbert scheme
$\alen$ also admits the same quiver variety description as
$\hquiver$ but with a different stability condition. It follows
from Nakajima's results on quiver varieties (Corollary 4.2,
\cite{Na1}) that $\hquiver$ is diffeomorphic to $\alen$.

When $n=1$, a theorem of Ito-Nakamura \cite{INr1} says that
$\varphi$ is an isomorphism. When $\G$ is trivial, $\phi:
\C^{2n}//S_n \rightarrow (\C^2)^{[n]}$ being an isomorphism is a
very difficult theorem due to Haiman \cite{Hai2}. Our construction
above of the morphism $\varphi$ was inspired by his construction.
Indeed Haiman's theorem is equivalent to the validity of the
Garsia-Haiman $n!$ conjecture, which in turn implies the Macdonald
positivity conjecture \cite{Hai1, Mac}. Our work provides a
natural geometric framework which may allow a generalization (even
for $\G$ cyclic) of these combinatorial results. This will be
pursued elsewhere.
\subsection{Equivalence of derived categories}

This subsection may be regarded as a continuation of
Subsect.~\ref{subsec_goodres}. I thank A.~King for stimulating
discussions. In this subsection, we assume that $Y$ is a
(quasi-)projective surface acted by a finite group $\G$ with
isolated singularities. Let $\tau: X \rightarrow Y/\G$ be
thecrepant resolution. Let $\tau_n: \Xn \rightarrow Y^n/\Gn$ be
the induced resolution of singularities.

Recall the Hilbert scheme of regular orbits $X^n//S_n$ was defined
in the previous subsection. Denote by $\mathcal Z$ the universal
closed subscheme ${\mathcal Z}\subset (X^n//S_n) \times X^n$ and
denote by $p, q$ the projections to $X^n$ and $X^n//S_n$. We have
the commutative diagram

$$
 \CD \mathcal Z @>{q}>> X^n \\ @VV{p}V @VVV \\ X^n//S_n @>>>
 X^n /S_n \endCD
$$
Haiman's theorem \cite{Hai2} allows us to identify $X^n//S_n \cong
X^{[n]}$, which we will use below implicitly. We denote by
$D_{\Gn} (Y^n)$ the bounded derived category of $\Gn$-equivariant
coherent sheaves on $Y^{n}$, and denote by $D(\Xn)$ the bounded
derived category of coherent sheaves on $\Xn$. Define two functors
$\Phi : D(\Xn) \rightarrow D_{S_n} (X^n)$ and $\Psi : D_{S_n}
(X^n) \rightarrow D(\Xn) $ by
\begin{eqnarray*}
 \Phi (-)  & =& Rp_* ({\mathcal O}_{\mathcal Z} \bigotimes q^*(-))\\
 \Psi (-)  & =& (Rq_* RHom({\mathcal O}_{\mathcal Z}, p^* (-)))^{S_n}.
\end{eqnarray*}
By applying a remarkable result of Bridgeland, King and Reid
\cite{BKR} together with Haiman's isomorphism $X^n//S_n \cong
X^{[n]}$, we see that $\Phi$ is an equivalence of categories and
$\Psi$ is its adjoint functor.

In a similar (and easier) way \cite{KV, BKR} (also cf.
\cite{GSV}), we also have an equivalence of categories $\phi: D(X)
\rightarrow D_\G(Y)$. This induces an equivalence of categories
$\phi^n: D(X^n) \rightarrow D_{\G^n}(Y^n)$ by passing to the
$n$-th Cartesian product. If we restrict ourselves to the
$S_n$-equivariant sheaves, we obtain from the equivalence $\phi^n$
an equivalence of categories between $D_{S_n}(X^n)$ and
$D_{\Gn}(Y^n)$ (note that a sheaf on $Y^n$ which is equivariant
with respect to $\G^n$ and $S_n$ is precisely a sheaf which is
equivariant with respect to the wreath product $\Gn =\G^n \rtimes
S_n$).

Combining with the equivalence $\Phi$, we have established the
following.

\begin{theorem}
Assume that $Y$ is a (quasi-)projective surface acted by a finite
group $\G$ with isolated singularities. Let $\tau: X \rightarrow
Y/\G$ be the crepant resolution. There exists an equivalence of
categories between $D_{\Gn}(Y^n)$ and $D(\Xn)$ associated to the
resolution of singularities $\tau_n: \Xn \rightarrow Y^n/\Gn$.
\end{theorem}

In other words, we can say that if $\tau: X \rightarrow Y/\G$ is
good in the sense of equivalence of derived categories then
$\tau_n: \Xn\rightarrow Y^n/\Gn$ is also good. This settles some
questions and confusing points in \cite{Wa2}. In the case when
$\G$ is trivial and $X$ equals $Y$, it reduces to the equivalence
between $D(\Xn)$ and $D_{S_n} (X^n)$.
\subsection{Cup product on $\Xn$ vs convolution product on $\Gn$}

So far we have seen there is a surprising analog (cf. Eqns.
(\ref{eq_bound}) and (\ref{eq_conv})) between the convolution
product $K_i(c,n)$ (at least when $i =0, 1$) on the representation
ring $R(\Gn)$ and the cup product with $B_i(\g, n)$ in the
cohomology ring $\Hn$. In the case when $\G$ is a finite subgroup
of $SL_2(\C)$, and $X$ is the minimal resolution $\ale$ of the
simple singularity $\C^2/\G$, we have expected a precise
connection between these two products which will generalize the
isomorphism of vector spaces between $H^*(\alen)$ and $R(\Gn)$.

In the case when $\G$ is trivial, Lehn and Sorger \cite{LS} indeed
found a precise and canonical connection between them. This result
has been also obtained by Vasserot \cite{Vas} by a different
technique. Let us associate a degree $d = n -\ell (\lambda)$ to a
permutation $\sigma$ of cycle type $\lambda$. This defines a
grading on $\C [S_n]$. The convolution product does not preserve
this grading but preserves the filtration on $R(S_n)$
\[
 F^0\C [S_n] \subset F^1 \C [S_n]  \subset \cdots \subset
 F^{n-1} \C [S_n] =\C [S_n],
\]
where $F^d \C [S_n]$ is the space spanned by permutations of
degree $\le d$. We call the induced product on $\C[S_n] =
\oplus_{d=0}^{n-1} F^d\C[S_n]/F^{d-1}\C[S_n]$ the {\em filtered
convolution product}, and denote it by $\cup$. In other words,
given two permutations $\sigma$ and $\pi$, we see that $\sigma
\cup \pi$ coincides with the usual convolution product $ \sigma
* \pi$ when $\deg (\sigma) +\deg (\pi) = \deg(\sigma \pi)$ and is $0$
otherwise.

We may restrict the filtered convolution to
$R(S_n) \subset \C[S_n]$. Recall a linear basis of $H^*(
(\C^2)^{[n]})$ is given by $\mathfrak q_{\lambda} =\prod_{r \ge 1}
\mathfrak q_r^{m_r} \vac$ while a linear basis of $R(S_n)$ is
given by $p_{\lambda} =\prod_{r \ge 1} p_r^{m_r} \vac$, where
$\lambda =(1^{m_1}2^{m_2}\ldots)$ runs over all partitions of $n$.
Here we also use $\vac$ to denote $1 \in R(\G_0) \cong \C$.

\begin{theorem} \cite{LS, Vas}
The linear map from cohomology ring $H^*( (\C^2)^{[n]})$ to the
filtered representation ring $R(S_n)$ by sending $\mathfrak
q_{\lambda}$ to $p_{\lambda}$ is a graded ring isomorphism.
\end{theorem}

Note that the formula (\ref{eq_bound}) holds for a projective
surface $X$. For the affine plane which is quasi-projective,
certain degeneracy occurs and the modified formula becomes
\cite{Lehn}

\begin{eqnarray} \label{eq_modify}
\mathfrak d =-\frac12 \sum_{n, m> 0}
   nm \mathfrak q_{n+m}\frac{\partial}{\partial \mathfrak q_{n}}
   \frac{\partial}{\partial \mathfrak q_{m}}.
\end{eqnarray}

Let us briefly comment on how Lehn and Sorger proved their
theorem. If one examines the proof of the equation
(\ref{eq_conv}), one sees that the operator acting on $R(S_n)$
which corresponds to $\mathfrak d$ given in (\ref{eq_modify}) is
the {\rm filtered} convolution product with the conjugacy class of
cycle type $(1^{n-2}2)$. With the help of the two parallel
formulas (\ref{eq_class}) and (\ref{eq_sign}), one can establish
that the cup product with $c(p_{1*}\mathcal O_{\mathcal Z_n} )$ in
$H^*( (\C^2)^{[n]})$ corresponds exactly to the filtered
convolution product with $\varepsilon_n$ in $R(S_n)$. These two
statements are mainly responsible for establishing the theorem
above. Indeed, Lehn-Sorger managed to go further to establish a
stronger form of the above theorem which claims the graded ring
isomorphism is valid between $H^*((\C^2)^{[n]}, \Z)$ and
$R_\Z(S_n)$.

In \cite{CR}, Chen and Ruan introduced the orbifold cohomology
ring for a general orbifold which is not necessarily a global
quotient by a finite group. Ruan observed that the filtered
convolution product on $R(S_n)$ above coincides with the orbifold
cup product on the orbifold cohomology ring of the orbifold
$\C^{2n}/S_n$, cf. \cite{Ruan}. Therefore by combining with the
theorem above, the cohomology ring on the Hilbert scheme
$(\C^2)^{[n]}$ is isomorphic to the orbifold cohomology ring of
the symmetric product $\C^{2n}/S_n$. This raises the question how
to relate the cohomology ring on the Hilbert schemes to the
orbifold cohomology ring on the symmetric product or more
generally on the wreath product orbifolds.
\section{Appendix: A quiver variety description of $\hquiver$}
\label{append_quiver}

We recall (cf. \cite{Na3}) that the Hilbert scheme $(\C^2)^{[K]}$
of $K$ points in $\C^2$ admits a description in terms of a quiver
consisting of one vertex and one arrow starting and ending at that
vertex. More explicitly, we denote
\begin{eqnarray*}
 \widetilde{H}(K) =
\left\{
\begin{array}{lcl}
              &|& i) [B_1, B_2 ] +ij =0  \\
(B_1, B_2, i, j) &|& ii) \mbox{ there exists no proper subspace
}\\
              &|&  S \subset
              \C^K \mbox{ such that } B_{\alpha} (S) \subset S \mbox{ and } \\
               &|& \mbox{ im } i \subset S \quad (\alpha =1,2)
\end{array}
              \right\}  ,
\end{eqnarray*}
where $B_1, B_2 \in \End (\C^K ), i \in \Hom (\C, \C^K ), j \in
\Hom (\C^K, \C )$. Then we have an isomorphism
\begin{eqnarray}  \label{eq_quot}
(\C^2)^{[K]} \cong \widetilde{H}(K)/ GL_K (\C ) ,
\end{eqnarray}
where the action of $GL_K (\C )$ on $\widetilde{H}(K)$ is given by
\begin{eqnarray*}
  g \cdot (B_1, B_2, i, j) = (g B_1 g^{ -1}, g B_2 g^{ -1}, g i, jg^{-1}).
\end{eqnarray*}
It is also often convenient to regard $(B_1, B_2)$ to be in $\Hom
(\C^K , \C^2 \otimes \C^K )$. We remark that one may drop $j$ in
the above formulation because one can show by using the stability
condition that $j =0$ (cf. \cite{Na3}, Proposition~2.8).

The bijection in (\ref{eq_quot}) is given as follows. For $I \in
(\C^2)^{[K]}$, i.e. an ideal in $\C [x, y]$ of colength $K$, the
multiplication by $x, y$ induces endomorphisms $B_1, B_2$ on the
$K$-dimensional quotient $\C [x, y] /I$, and the homomorphism $i
\in \Hom (\C, \C^K)$ is given by letting $i (1) =1 \mbox{ mod }
I$. Conversely, given $(B_1, B_2, i)$, we define a homomorphism
$\C [x, y] \rightarrow \C^K$ by $f \mapsto f(B_1, B_2) i (1)$. The
stability condition is equivalent to the homomorphism $\C [x, y]
\rightarrow \C^K$ being surjective, which implies that the kernel
$I$ of this homomorphism is an ideal of $\C [x, y]$ of length $K$.
One easily checks that the two maps are inverse to each other.

Set $K =nN$, where $N$ is the order of $\G$. We may identify
$\C^K$ with the direct sum of $n$ copies of the regular
representation $R$ of $\G$, $\C^2$ with the defining
representation $Q$ of $\G$ by the embedding $\G \subset SL_2
(\C)$, and $\C$ with the trivial representation $V_{\g_0}$ of
$\G$. Denote by
\begin{eqnarray*}
M (n) & =& \Hom (R^n , Q \otimes R^n ) \bigoplus \Hom (V_{\g_0},
R^n)\bigoplus \Hom (R^n, V_{\g_0}).
\end{eqnarray*}
By definition $\widetilde{H} (nN) \subset M(n)$. Let $GL_{\G} (R)$
be the group of $\G$-equivariant automorphisms of $R$. Then the
group $G \equiv GL_{\G} (R^n)$ acts on the $\G$-invariant subspace
$M(n)^{\G} $. We give the following description of $\hquiver$ as a
quiver variety. This result is certainly known to Nakajima and it
seems to be observed by others as well. A proof of it is given in
\cite{Wa2}.
\begin{theorem}  \label{th_quiv}
The variety $\hquiver $ admits the following description:
\begin{eqnarray*}
  \hquiver \cong (\widetilde{H} ( nN) \cap M(n)^{\G} ) / GL_{\G} (R^n).
\end{eqnarray*}
\end{theorem}

Consider the $\G$-module decomposition $Q \otimes V_{\g_i}=
\bigoplus_j a_{ij} V_{\g_j}$, where $a_{ij} \in \Z_+$, and
$V_{\g_i} $ $(i =0, \ldots, r)$ are irreducible representations
corresponding to the characters $\g_i$ of $\G$. Set $\dim V_{\g_i}
= n_i$. Then
\begin{eqnarray}
 M(n)^{\G}
 & =& \Hom_{\G} (R^n ,Q \otimes R^n )
  \bigoplus \Hom_{\G} (V_{\g_0}, R^n)
  \bigoplus \Hom_{\G} (R^n, V_{\g_0})
 \label{eq_invt} \\
 &\cong& \Hom_{\G} (\sum_i \C^{n n_i} \otimes V_{\g_i},
  \C^2 \otimes \sum_i \C^{n n_i} \otimes V_{\g_i} )\nonumber \\
  && \bigoplus \Hom_{\G} (V_{\g_0}, R^n)
   \bigoplus \Hom_{\G} (R^n, V_{\g_0}) \nonumber \\
 &\cong& \sum_{i j} a_{ij} \Hom ( \C^{n n_i} ,\C^{n n_j})
  \bigoplus \Hom (\C, \C^n )
   \bigoplus \Hom (\C^n, \C).  \nonumber
\end{eqnarray}
where $\Hom_{\G}$ stands for the $\G$-equivariant homomorphisms.
In the language of quiver varieties as formulated by Nakajima
\cite{Na1, Na4}, the above description of $\hquiver$ identifies
$\hquiver$ with a quiver variety associated to the following data:
the graph consists of the same vertices and edges as the McKay
quiver which is an affine Dynkin diagram associated to a finite
subgroup $\G$ of $SL_2 (\C)$; the vector space $V_i$ associated to
the vertex $i$ is isomorphic to the direct sum of $n$ copies of
the $i$-th irreducible representation $V_{\g_i}$; the vector space
$W_i =0$ for nonzero $i$ and $W_0 =\C$.

The variety $\hquiver$ is connected since a quiver variety is
always connected \cite{Na1, Na4}. According to Nakajima
(unpublished), the Hilbert scheme $\alen$ admits a quiver variety
description in terms of the same quiver data as specified in the
above paragraph but with a different stability condition. It
follows by Corollary 4.2 of \cite{Na1} that $\alen$ and $\hquiver$
is diffeomorphic.

\pagebreak
 \vspace{1cm}
\begin{table}
 \begin{center}
\caption{A DICTIONARY} \label{Tab_dict}
\begin{tabular}{|| l | l || }
\hline
 Hilbert Scheme $\Xn$ Picture &  Wreath Product $\Gn$ Picture\\ \hline\hline
   length-$n$ schemes supported at a point
                              &  $n$-cycle   \\ \hline
 cohomology group $\Hn$       &  Grothendieck group $R(\Gn)$   \\ \hline
 $\Hx \cong$ Heisenberg Fock space
                              &  $\RG \cong$ Heisenberg Fock space  \\ \hline
 lattice $H^*(X,\Z)/{tor}$
                              &  lattice $R_\Z(\G)$  \\ \hline
   $ (\C^2)^{[n]} $           &  $S_n$ \\ \hline
   $ \alen$ for $\G \subset SL_2(\C)$
                              &  $\Gn$ for $\G \subset SL_2(\C)$ \\ \hline
 cup product ($X$ closed)     &  convolution product   \\ \hline
 cup product ($X$ non-closed) &  (filtered) convolution product   \\ \hline
 correspondence varieties     &  induction/restriction functors etc  \\ \hline
 boundary of $\Xn$            & conjugacy class of type $(1^{n-2}2)$    \\ \hline
 cohomology class $B_i(\alpha,n)$
                              &  conjugacy class $K_i(c, n)$   \\ \hline
 total Chern class $c(L^{[n]})$
                              &  character $\varepsilon_n(\g)$  \\ \hline
   $?$                        &  character $\eta_n(\g)$   \\ \hline
 cohomology class $G_i(\alpha,n)$
  &  ?   \\ \hline

\end{tabular}
 \end{center}
\end{table}

\end{document}